\documentclass{amsart}
\usepackage{amssymb}
\newcommand{\0}{{\mathbf  0}}
\newcommand{\1}{{\mathbf 1}}
\newcommand{\rationals}{\mathbb Q}
\newcommand{\N}{{\mathcal N}}

\newcommand{\add}{\text{\sf{add}}}
\newcommand{\reals}{\mathbb R}
\newcommand{\Proof}{{\sc Proof} \hspace{0.2in}}
\newcommand{\rest}{\restriction}
\newcommand{\<}{\langle}
\renewcommand{\>}{\rangle}
\newcommand{\lft}[2]{\mathopen\ifcase#1{}\oo\or
                        \big#2\or\Big#2\else\oo\fi} 
\newcommand{\rgt}[2]{\mathclose\ifcase#1{}\oo\or
                        \big#2\or\Big#2\else\oo\fi} 
\newcommand{\QED}{\hspace{0.1in} \square \vspace{0.1in}}
\newcommand{\dom}{{\text{\sf {dom}}}}

\theoremstyle{plain}
\newtheorem{theorem}{Theorem}[section]
\theoremstyle{plain}
\newtheorem{lemma}[theorem]{Lemma}
\newtheorem{definition}[theorem]{Definition}

\begin{document}

\title{Not every $\gamma$-set is strongly meager}
\author{Tomek Bartoszy\'{n}ski}
\address{Department of Mathematics\\
   Boise State University\\
   Boise, Idaho 83725}
\email{tomek@@math.idbsu.edu}
\author{Ireneusz Rec\l aw}
\address{Institute of Mathematics\\ 
University  Gdansk \\
Wita Stwosza 57\\
 80-952 Gdansk, Poland}
 \email{matir@@halina.univ.gda.pl}
\thanks{First author partially supported by  SBOE grant \#95--041 and
  the second author supported by the research grant
BW UG 5100-5-0148-4}
\subjclass{04A15 03E50}

\maketitle
\begin{abstract}
We present two constructions of $\gamma$-sets which are large in sense
of category.
\end{abstract}
\section{Introduction}
Most of the notation used in this paper is standard. By reals we mean
the space $2^\omega$ with the operation $+ \mod 2$. By rationals we
mean the canonical dense subset of $2^\omega$, 
$$\rationals=\{x \in 2^\omega : \forall^\infty n \ x(n)=0\}.$$
Finally, let $\N$ denote the ideal of measure zero subsets of
$2^\omega$, with respect to the standard product measure in this
space.
$\0$ and $\1$ denote constant functions equal to $0$ and $1$
respectively. 
 
  Let us recall  the following definitions:
\begin{definition}
\begin{enumerate}
\item     A family ${\mathcal J} \subseteq P(X)$ is an $\omega$-cover of
     $X$
 if for  every  finite 
set $F \subseteq X$ there exists $B \in {\mathcal J}$ such that $F \subseteq B$,
     \item A topological space $X$ is a $\gamma$-set if for  every
       ${\mathcal J}$,   open 
$\omega$-cover of $X$, there exists a family $\{D_n : n \in \omega\}
\subseteq {\mathcal J}$
 such that 
$X \subseteq \bigcup_m \bigcap_{n>m} D_n$,
\item A set $X \subseteq \reals$ is strongly meager if for every null set
$G\subset \reals$, $X + G \ne \reals$.
\end{enumerate}
\end{definition}

Galvin and Miller in \cite{GalMil84Gam} constructed a $\gamma$-set of
reals of size continuum  
under Martin's Axiom. They showed that for every $\gamma$-set $X \subseteq 
\reals$ and
 every meager set $F$, $X+F$ is meager. They asked 
whether the same is true
for null sets. We give negative answer to this question.

Note that $\gamma$-sets are always meager, that is, if $X$ is a
$\gamma$-set  then $X \cap P$ is meager in $P$ for every perfect set
$P \subseteq 2^\omega$ (\cite{GN}).

 Since every $\gamma$-set
is a strong measure zero set (see \cite{GN}), by a result of Laver it
is consistent  
that every $\gamma$-set is countable so a sum with a null set is null.

On the other hand, for an ideal ${\mathcal I}$, Pawlikowski defined the cardinal 
coefficients $\add_t({\mathcal I})=\min\{|X|:
\forall {F \in {\mathcal I}} \ X + F \not\in {\mathcal I}\}$.
He showed in \cite{Paw}, 
that $\add(\N)=\min\{\add_t(\N),{\mathfrak b}\}$. It is consistent
that $\add(\N)<{\mathfrak p}$. Then $\add_t(\N)=\add(\N)<{\mathfrak p}$. 
So every subset
of size $\add_t(\N)$ is a $\gamma$-set, since all sets of size $<
{\mathfrak p}$ are $\gamma$-sets (\cite{Mil84Spe}), but there is a set $X$ of  size
$\add_t(\N)$ and 
a null set $G$ with $X+G \not\in \N$.

\section{Main results}

In this section we present two constructions of $\gamma$-sets.

\begin{theorem}\label{thm1}
Assume ${\mathfrak p} = 2^{\boldsymbol\aleph_0}$. Then there is 
a $\gamma$-set $X \subseteq 2^\omega$ which is not strongly meager.
\end{theorem}
\Proof
For $n \in \omega$ let $k_n=\sum_{i=0}^n 2^i$. 
Let $A_n=\{x \in 2^\omega: x \rest [k_n,k_{n+1})= \0 \}$ for $n \in \omega$.
Note that each set $A_n$ is clopen and has measure 
 measure $2^{-(n+1)}$. So $G=\bigcap_m
\bigcup_{n>m} A_n$ 
is null.

We will construct a set $X'=\{x_\alpha:\alpha < 2^\omega\}$ such 
that 
\begin{enumerate}
\item $\forall \alpha \ x_\alpha \in 2^\omega$,
\item $\forall \alpha \ \forall^\infty n \ x_{\alpha+1}(n) \leq x_\alpha(n)$,
\item $X=X'\cup \rationals$
is a $\gamma$-set,  
\item $\forall z \in 2^\omega \  \exists \alpha \ 
\exists^\infty n\  z\rest [k_n,k_{n+1})= x_\alpha \rest [k_n,k_{n+1})$,
\item $\forall \alpha \ \exists^\infty n \ x_\alpha \rest
  [k_n,k_{n+1})= \1$.
\end{enumerate}

Note that if $z$ and $x_\alpha$ are like in (4) then $z+x_\alpha \in
G$. Therefore, condition (4) implies that $X+G =2^\omega$.

Let ${\mathcal S}=\{z \in 2^\omega: \exists^\infty n \ 
z\rest [k_n,k_{n+1}) =\1\}$.
For $x \in 2^\omega$ let  $[x]^\star=\{z \in {\mathcal S}: \forall^\infty n \
x(n) \geq z(n) \}$.

\begin{lemma}\label{lemma1}
 Let ${\mathcal J}$ be open $\omega$-cover of $\rationals$
and $x \in {\mathcal S}$.
Then there is a sequence $D_n \in {\mathcal J}$ and $y \in [x]^\star$ such that
$\rationals \cup [y]^\star \subseteq \bigcup_m \bigcap_{n>m} D_n$.
\end{lemma}
\Proof
 Using the fact that ${\mathcal J}$ is open $\omega$-cover of $\rationals$
 we can find a sequence $\<l_n: n \in \omega\>$ of natural numbers and
 a sequence $\<D'_n: n \in \omega\>$ of elements of ${\mathcal J}$ such
 that for every $n$,
$$\forall z\in 2^\omega  \ \lft2(z \rest [n,l_n) = \0 \rightarrow z \in D'_n\rgt2).$$
Without loss of generality we can assume that there exists a set $Z
\subseteq \omega$ such that 
$$x \rest [k_n, k_{n+1}) = \left\{\begin{array}{ll}
\1 & \text{if $n \in Z$}\\
\0 & \text{if $n \not\in Z$}\end{array}\right. .$$
Choose $Y \subseteq \omega$ and $Z' \subseteq Z$ such that 
$$\bigcup_{n \in Z'} [k_n, k_{n+1}) \cap \bigcup_{n \in Y} [n, l_n) =
\emptyset.$$
Define
$$y \rest [k_n, k_{n+1}) = \left\{\begin{array}{ll}
\1 & \text{if $n \in Z'$}\\
\0 & \text{if $n \not\in Z'$}\end{array}\right. .$$
It is clear that $y \in [x]^\star$. Suppose that $z \in [y]^\star$.
Note that
$$\forall^\infty n \in Y \ z \rest [n,l_n) = \0,$$
which means that $z \in D'_n$ for all except finitely many $n \in Y$. 
Thus, in order to finish the proof it is enough to define
$D_n = D'_{y(n)}$, where $y(n)$ is the $n$-th element of $Y$.~$\QED$

\begin{lemma}\label{lemma2}
 Suppose that $\{x_\alpha: \alpha<\kappa< {\mathfrak p}\} \subseteq
 2^\omega$  is a sequence such that $x_\alpha \in [x_\beta]^\star$ for
 $\alpha>\beta$. Then
$\bigcap_{\alpha<\kappa} [x_\alpha]^\star \neq \emptyset$.
\end{lemma}
\Proof  Define $Y_\alpha =\left\{n: x_\alpha \rest [k_n,k_{n+1})= \1 \right\}$.
Then $Y_\alpha  \subseteq ^\star Y_\beta$ if $\alpha  > \beta$. Since
$\kappa<{\mathfrak p}$, there
is $Y$ such that $Y \subseteq^\star Y_\alpha$ for each $\alpha  <
\kappa$.
Define 
$$x \rest [k_n, k_{n+1}) = \left\{\begin{array}{ll}
\1 & \text{if $n \in Y$}\\
\0 & \text{if $n \not\in Y$}\end{array}\right. .$$
It is clear that $x \in \bigcap_{\alpha<\kappa} [x_\alpha]^\star$.~$\QED$

Let $\{{\mathcal J}_\alpha: \alpha < 2^\omega \}$ be an enumeration of 
all $\omega$-covers of 
$\rationals$, and let $\{z_\alpha: \alpha < 2^\omega\}$ be enumeration
of all elements of $2^\omega$.

Assume that the set $X_\alpha=\{x_\beta:\beta < \alpha \}$ has been already
constructed. 
Assume that ${\mathcal J}_\alpha $ is $\omega$-cover of $X_\alpha \cup
\rationals$. Since $|\alpha|<{\mathfrak p}$ and all sets of size $< {\mathfrak p}$
are $\gamma$-sets,
we can choose ${\mathcal J}'_\alpha=\{U_n: n \in \omega\}$, $\omega$-subcover of 
this such that
$$X_\alpha \cup \rationals \subseteq \bigcup_{m} \bigcap_{n>m} U_n.$$

If $\alpha$ is limit 
then apply 
\ref{lemma2} to get   a real $x_\alpha' \in \bigcap_{\beta<\alpha}
[x_\beta]^\star$. If $\alpha $ is
not limit let $x_\alpha' = x_{\alpha-1}$.

Next apply \ref{lemma1} to $x_\alpha'$ and ${\mathcal J}'_\alpha$ and
$x_\alpha'$ to get a real $y_\alpha$. 
Finally let $x_\alpha$ be such that 
\begin{enumerate}
\item $\exists^\infty n\  z_\alpha \rest [k_n,k_{n+1})= x_\alpha \rest [k_n,k_{n+1})$,
\item $ \exists^\infty n \ x_\alpha \rest
  [k_n,k_{n+1})= \1$,
\item $\forall^\infty n \ x_\alpha(n) \leq y_\alpha(n)$.
\end{enumerate}
This finishes the construction and the proof of the theorem.~$\QED$

The set constructed above is
a $\gamma$-set but it contains a subset which is not a $\gamma$-set. In
the next theorem we will show how to build a set which is a hereditarily
$\gamma$-set and is not strongly meager. 
The construction is a slight modification of Todorcevic's construction
of a hereditarily $\gamma$-set from \cite{GalMil84Gam}.
\begin{theorem}\label{thm2}
Assume $\diamond$. Then there exists a hereditarily $\gamma$-set which
is not strongly meager.
\end{theorem}
\Proof
For a tree $p \subseteq 2^{<\omega}$ let  $[p]$ be 
the set of branches of $p$. Similarly, for a finite set $U
\subseteq 2^{<\omega}$ let 
$$[U]= \{x \in 2^\omega : \exists s \in U \
x \rest \dom(s)=s\}.$$

Let $\{k_n: n \in \omega\}$ and $G$ be the sequence and the set
defined at the beginning of the proof of \ref{thm1}.
Let ${\mathcal P}$ be the collection of all perfect trees $p$ such that 
there exists a sequence $\{U_n: n \in \omega\}$ such that
\begin{enumerate}
\item $U_n \subseteq 2^{[k_n, k_{n+1})}$ for all $n$,
\item $\exists^\infty n \ U_n = 2^{[k_n, k_{n+1})}$,
\item $[p]= \bigcap_{n \in \omega} [U_n]$.
\end{enumerate}
By induction on levels we will build an Aronszajn tree consisting of
elements of ${\mathcal P}$ 
ordered by inclusion. The $\gamma$-set we are looking for will be a
selector from the elements of this tree.

For perfect trees $p,q$ and a set $R \subseteq 2^n \cap q$ define
$$p \leq_R q \Rightarrow p \cap 2^n =R \ \&\ p \subseteq q.$$
If $R=q \cap 2^n$ we write $p \leq_n q$ instead of $p \leq_{q \cap
  2^n} q$.

Let $\{z_\alpha: \alpha<\omega_1\}$ be enumeration of $2^\omega$.

We will build by induction a partial ordering $\prec$ on $\omega_1$, 
$\{p_\alpha: \alpha \in \omega_1\}$ and $\{x_\alpha:\alpha \in \omega_1\}$ such that 
\begin{enumerate}
\item $T=(\omega_1, \prec)$ is an Aronszajn tree and for limit $\alpha$,
  $\bigcup_{\beta<\alpha} T_\beta=\alpha$,
\item $p_\alpha \in {\mathcal P}$ for all $\alpha$,
\item $\forall \alpha, \beta \ \lft2(\alpha \prec \beta \iff p_\alpha
  \subseteq p_\beta\rgt2)$, 
\item $x_\alpha \in [p_\alpha]$ for all $\alpha$,
\item $\exists^\infty n \ x_\alpha \rest [k_n, k_{n+1})=z_\alpha \rest
  [k_n, k_{n+1})$, 
\item if ${\mathcal J}$ is an ``appropriate'' 
$\omega$-cover then there exists $\alpha$
  and a sequence $\<D_n:n \in \omega\>$ of elements of ${\mathcal J}$
  such that  for all $p \in T_\alpha$
$[p] \subseteq \bigcup_{m} \bigcap_{n>m} D_n$,
\item if $\beta>\alpha$ and $q \in T_\alpha$, $R \subseteq q \cap 2^m$
 then there exists $p \in T_\beta$  such that $p \leq_R q$.
\end{enumerate}
Note that the condition $(7)$ guarantees that the construction will
not terminate after countably many steps.
Condition $(6)$ is rather vague, but with the right interpretation of
the word ``appropriate'', together with the condition $(4)$ 
it will guarantee that the set
$X=\{x_\alpha: \alpha < \omega_1\}$  is a hereditary $\gamma$-set.
Finally $(5)$ yields that $X+G=2^\omega$.

For a set $Y \subseteq 2^\omega$ and a perfect tree $p$ let 
$p(Y)$ be the tree representing the closure of $[p] \cap Y$.

We  use $\diamond$ to construct oracle sequences $\{{\mathcal J}_\alpha,
X_\alpha, \<p^\beta: \beta<\alpha\>: \alpha <\omega_1\}$ such that 
for any $Y \subseteq X$ and an open $\omega$-cover of $Y$, ${\mathcal J}$, there
are stationary many  $\alpha$'s such that
\begin{enumerate}
\item ${\mathcal J}_\alpha={\mathcal J}$,
\item $X_\alpha = \{x_\beta: \beta<\alpha\} \cap Y$,
\item $p^\beta=p_\beta(Y)$ for $\beta<\alpha$.
\end{enumerate}
Note that these sequences are easy to obtain from an ordinary
$\diamond$ sequence by coding ${\mathcal J}$'s and $Y$'s by subsets of
$\omega_1$ (even $\omega$ in case of ${\mathcal J}$).

We will build the tree $T$ (or $\prec$) by induction on levels. If
$\alpha=\beta+1$ and $T_\beta$ is already constructed then $T_\alpha$
is any extension of $T_\beta$ satisfying the requirements.

Suppose that $\alpha$ is a limit ordinal. We look for a sequence
$\{D_n: n \in \omega\} \subseteq {\mathcal J}_\alpha$ such that 
\begin{enumerate}
\item[(i)] $\forall y \in X_\alpha \ \forall^\infty n \  y \in D_n$,
\item[(ii)] $\forall \beta <\alpha \ \forall m \ \forall R \subseteq p_\beta
  \cap 2^m \ \exists \delta \ \lft2(p_\delta \leq_R p_\beta \ \&\
  \forall^\infty n \ [p^\delta] \subseteq D_n\rgt2)$.
\end{enumerate}
If such a sequence $\{D_n:n \in \omega\}$ does not exist then
$T_\alpha$ is an arbitrary extension of $\bigcup_{\beta<\alpha}
T_\beta$. 

Otherwise we fix a sequence $\{D_n:n \in \omega\}$ satisfying  the above
conditions then for every $\beta<\alpha$ and every $R \subseteq
p_\beta \cap 2^m$ we build a chain $\{\delta_n: n \in \omega\}
\subseteq \alpha$ and $\{l_n:n \in \omega\}$ such that 
\begin{enumerate}
\item $\delta_n \prec \delta_{n+1}$ for all $n$,
\item $p_{\delta_{n+1}} \leq_{l_n} p_{\delta_n}$ for all $n$,
\item $\bigcap_{n
  \in \omega} p_{\delta_n} \in {\mathcal P}$,
\item $p_{\delta_0} \leq_R p_\beta$,
\item $\forall^\infty n \ [p^{\delta_0}] \subseteq D_n$.
\end{enumerate}
The branch $\{\delta_n:n \in \omega\}$ will be extended on level
$\alpha$ by, say, $\rho$ and the corresponding set is $p_\rho=\bigcap_{n
  \in \omega} p_{\delta_n}$. 
Note that in this way we extend only countably many branches.

Observe that condition $(3)$ will follow from $(2)$ if only the sequence
$l_n$ is increasing fast enough. Conditions $(4)$ and $(5)$ follow from
the condition (ii) above. This concludes the construction.
It remains to show that $X$ is a hereditary $\gamma$-set. 

Suppose that $Y \subseteq X$ and let ${\mathcal J}$ be an $\omega$-cover
of $Y$. Let $\alpha$ be a limit ordinal such that 
\begin{enumerate}
\item ${\mathcal J}_\alpha={\mathcal J}$,
\item $X_\alpha = \{x_\beta: \beta<\alpha\} \cap Y$,
\item $p^\beta=p_\beta(Y)$ for $\beta<\alpha$.
\end{enumerate}
To finish the proof it is enough to check that there exists a sequence
$\{D_n: n \in \omega\} \subseteq {\mathcal J}_\alpha$ satisfying
conditions (i) and (ii). In this case, according to the construction, 
$Y \subseteq \bigcup_m
\bigcap_{n>m} D_n$. 

Let $\{\<\beta_n, R_n\>: n \in \omega\}$ be enumeration of the set 
$$\{\<\beta, R\>: \beta<\alpha \ \&\ \exists m \ R \subseteq
p_\beta \cap 2^m\}.$$

It is enough to construct by induction a sequence $\{D_n: n \in
\omega\} \subseteq {\mathcal J}_\alpha$ such that 
\begin{enumerate}
\item[(iii)] $\forall i<n \ x_{\beta_i} \in D_n$,
\item[(iv)] $\forall i < n \ \exists \delta \ p_\delta \leq_{R_i}
  p_{\beta_i} \ \&\ [p^\delta] \subseteq D_n$.
\end{enumerate}

We describe how to construct  set $D_n$. For each $i<n$ and $s \in
R_i$ choose $x^i_s \in Y$ (if possible) such that $s \subseteq x^i_s
$.
Let $D_n$ be an element of ${\mathcal J}_\alpha$ such that 
$$\{x_{\beta_i}: i<n\} \cup \{x^i_s: i<n, s \in R_i\} \subseteq D_n.$$
Such $D_n$ exists since ${\mathcal J}_\alpha$ is an $\omega$-cover of $Y$.

We need to verify condition (iv). Fix $i <n$ and 
choose $m$ so large that $[x^i_s \rest m] \subseteq D_n$ for all $s
\in R_i$. If $x^i_s$ does not exist let $x_s$ be any element of
$[p_{\beta_i}]$ extending $s$.
Let $R=\{x^i_s \rest m: s \in R_i\} \cup \{x_s \rest m: s \in R_i\}$.
By inductive hypothesis there exists $\delta$ such that $p_\delta
\leq_{R} p_{\beta_i}$. 
Note that $[p^\delta] \subseteq \bigcup_{s \in R_i} [x^i_s \rest
m]$. Thus $[p^\delta] \subseteq D_n$, which
finishes the construction and the proof.~$\QED$

We can show the same theorems for the algebraic structure of the real line.
Let us take a standard function $f: 2^\omega  \to [0,1]$,
$f(x)=\sum _{i \in \omega} x(i)/2^i$. $f$ is continuous so $f(X)$ is a
$\gamma$-set, 
where $X$ is as in \ref{thm1} or \ref{thm2}. 
Let $H= \bigcap_m \bigcup_{n>m}f(A_n)$. 
It is easy to see that $[0,1] \subseteq f(X) + H$ . Let $G= H + \rationals $
then $f(X)  + G = \reals$.~$\QED$

\ifx\undefined\bysame
\newcommand{\bysame}{\leavevmode\hbox to3em{\hrulefill}\,}
\fi

\end{document}